\documentclass[12pt]{article}
\usepackage{amsmath,amssymb,amsthm,latexsym,amscd}

\title{Forcing of infinity \\ and algebras of distributions \\
of binary semi-isolating formulas \\ for strongly minimal
theories\footnote{{\em Mathematics Subject Classification.} 03C07,
03G15, 20N02, 08A02, 08A55.
\newline\indent \ \ \ The work is supported by RFBR (grant
12-01-00460-a).}}
\author{Sergey V.
Sudoplatov\footnote{sudoplat@math.nsc.ru}}
\date{June 05, 2013}
\begin{document}
\maketitle

\begin{abstract}
We apply a general approach for distributions of binary isolating
and semi-isolating formulas to the class of strongly minimal
theories.

{\bf Key words:} structure of binary semi-isolating formulas,
strongly minimal theory.
\end{abstract}

\medskip
Algebras and structures associated with isolating and
semi-isolating for\-mu\-l\-as of a theory are introduced in
\cite{ShS, Su122}. We apply this general approach for
distributions of formulas to the class of strongly minimal
theories \cite{BaLa}.

\medskip
Let $T$ be a theory, $R\subseteq S^1(\varnothing)$ be a nonempty
family, and $\nu(R)$ be a regular family  of labelling functions
for semi-isolating formulas forming a set $U$ of labels. Denote by
$U_{\rm fin}$\index{$U_{\rm fin}$} (respectively $U_{\rm
cofin}$)\index{$U_{\rm cofin}$} the set of labels $u$, each of
which, being in $\bigcup\limits_{p,q\in R}\rho_{\nu(p,q)}$,
$\unlhd$-dominates a (co)finite set of labels in
$\rho_{\nu(p,q)}$.

By the definition all almost deterministic labels belong to
$U_{\rm fin}$.

We say that a label $u\in \rho_{\nu(p,q)}$ {\em
forces}\index{Label!forcing an infinite set of solutions} an
infinite set of solutions for the formula $\theta_{p,u,q}(a,y)$,
$\models p(a)$, if for any theory $T$ with a family $R$ of
$1$-types, containing $p$ and $q$ and having a ${\rm
POSTC}_R$-structure, including all labels that are
$\unlhd$-dominated by $u$, the formula $\theta_{p,u,q}(a,y)$ has
infinitely many solutions.

By the definition a label $u\in \rho_{\nu(p,q)}$ forces an
infinite set of solutions if and only if for any $n\in\omega$ and
some $n$ elements $a_1,\ldots,a_n$ in the set of solutions of
$\theta_{p,u,q}(a,y)$ (in an arbitrary structure), where $\models
p(a)$, there exists new element $a_{n+1}$, for which
$\models\theta_{p,u,q}(a,a_{n+1})$ and links between $a_{n+1}$ and
$a_1,\ldots,a_n$ are defined by some labels.

Clearly, almost deterministic labels do not force infinite sets of
solutions and any label $u$, $\unlhd$-dominating infinitely many
labels $v_i$, forces an infinite set of solutions. Moreover, each
label $u\wedge\neg v_{i_1}\wedge\ldots\wedge \neg v_{i_n}$ also
forces an infinite set of solutions. Another example with a label,
which forces an infinite set of solutions, is presented by theory
${\rm Th}(\langle\mathbb{Q},<\rangle)$, for which any non-zero
label (defining the strict order property) corresponds to formulas
having only infinitely many solutions. An infinite set of
solutions can be forced by labels in $U_{\rm fin}$ for formulas in
stable theories. Such an example is produced by any label
corresponding to a special element of an infinite group for an
everywhere finitely defined polygonometry \cite{SuGP}.

\medskip
{\bf Definition} (J.~T.~Baldwin\index{Baldwin J. T.},
A.~H.~Lachlan\index{Lachlan A. H.} \cite{BaLa}). A theory $T$ is
called {\em strongly minimal}\index{Theory!strongly minimal} if
for any formula $\varphi(x,\bar{a})$ of language obtained by
adding parameters of $\bar{a}$ (in some model ${\cal M}\models T$)
to the language of $T$, either $\varphi(x,\bar{a})$, or
$\neg\varphi(x,\bar{a})$ has finitely many solutions.

\medskip
An example of strongly minimal theory with the forcing of infinite
set is represented by structure $\langle M;s\rangle$ with {\em
successor function}\index{Function!successor} $s$ (having exactly
one preimage for any element, and do not having cycles). Since
${\rm Th}(\langle M;s\rangle)$ has unique $1$-type, there is a
label $u\in U_{\rm cofin}$ for the semi-isolating formula
$(x\approx x)$. This label $\unlhd$-dominates infinitely many
labels $v$ corresponding to formulas $(y\approx s^n(x))$,
$n\in\mathbb Z$, and thus, $u$ forces an infinite set of solutions
for the formulas $\theta_u(a,y)$, $a\in M$.

\medskip
{\bf Theorem 1.} {\em For any strongly minimal theory $T$, the
family $R\rightleftharpoons S^1(\varnothing)$ of $1$-types, and a
regular family $\nu(R)$ of labelling functions for semi-isolating
formulas, the following conditions hold:

{\rm (a)} the ${\rm POSTC}_\mathcal{R}$-structure
$\mathfrak{M}_{\nu(R)}$ consists of labels belonging to $U_{\rm
fin}\cup U_{\rm cofin}$;

{\rm (b)} there is unique type $r_0\in R$ having infinitely many
realizations; in particular, any set $\rho_{\nu(p,q)}$ is finite,
where $p,q\in R$ and $q\ne r_0$, and all labels
$u\in\rho_{\nu(p,q)}$ are almost deterministic and belong to
$U_{\rm fin}$;

{\rm (c)} if $R$ is finite, i.~e., all types in $R$ are principal,
then all non-zero labels are positive and all labels $u$,
including zero, have complements $\bar{u}$, and for any pair of
labels $u,\bar{u}\in\rho_{\nu(p,r_0)}$, exactly one of them is
almost deterministic and, in particular, belongs to $U_{\rm fin}$,
and the other label marks a formula $\theta_{p,\dot,r_0}(a,y)$
with infinitely many solutions, where $\models p(a)$, and belongs
to $U_{\rm cofin}$;

{\rm (d)} if $R$ is infinite, i.~e., $r_0$ is unique non-principal
$1$-type, then all non-zero labels, linking realizations of $r_0$
or realizations of types in $R\setminus\{r_0\}$, are positive, and
labels, linking realizations of $r_0$ with realizations of types
in $R\setminus\{r_0\}$, are negative; in this case, if a label $u$
belongs to $\rho_{\nu(p,r_0)}$ then $u$ is positive or zero,
almost deterministic, does not have complements and $p=r_0$,
moreover, $U_{\rm fin}=U_{\rm cofin}$ if $\rho_{\nu(r_0)}$ is
finite, and $U_{\rm cofin}=\varnothing$ if $\rho_{\nu(r_0)}$ is
infinite;

{\rm (e)} only labels in $\rho_{\nu(p,r_0)}$ with the principal
type $r_0$ can force infinitely many solutions.}

\medskip
{\em Proof.} By the definition of strongly minimal theory, each
formula $\varphi(a,y)$, where $\models p(a)$, (in particular,
witnessing the semi-isolation for a type $q(y)$) has a finite or a
cofinite set $\varphi(a,{\cal M})$ of solutions. For the finite
set $\varphi(a,{\cal M})$, the label $u\in\rho_{\nu(p,q)}$,
marking the formula $\varphi(x,y)$ with $\varphi(a,y)\vdash q(y)$
and $\models p(a)$ can $\unlhd$-dominate only finitely many labels
in $\rho_{\nu(p,q)}$. If $\varphi(a,{\cal M})$ is cofinite, then
the label $u\in\rho_{\nu(p,q)}$ for the formula $\varphi(x,y)$
$\unlhd$-dominates all labels in $\rho_{\nu(p,q)}$ except for a
finitely many $u_1,\ldots,u_k$, and $u$ has the complement
$\bar{u}$, which is obtained from $u$ by disjunctive attachment of
labels $u_i$. Thus the condition (a) holds: all labels belong to
$U_{\rm fin}\cup U_{\rm cofin}$.

Since $T$ is strongly minimal we also have that there is unique
type $r_0\in R$, principal or non-principal, with infinitely many
realizations: if there are finitely many $1$-types it is implied
by the property that models are infinite and there are no two
principal formulas with infinitely many realizations, and if there
are infinitely many $1$-types, then the non-principal type
$r_0(x)$, having infinitely many realizations by Compactness, is
isolated by the set of all formulas $\neg\varphi(x)$, where
$\varphi(x)$ are principal formulas and none of these formulas can
not have infinitely many solutions.

Since the type $r_0$ with infinitely many realizations is unique,
then any set $\rho_{\nu(p,q)}$ is finite, where $p,q\in R$ and
$q\ne r_0$. Here all labels $u\in\rho_{\nu(p,q)}$ are almost
deterministic and belong to $U_{\rm fin}$. Thus, we have the
condition (b).

If $r_0$ is isolated then all $1$-types are isolated and by
Proposition 1.1 \cite{Su122} all non-zero labels are positive.
Since each isolating formula $\varphi(x)$ has a label, all labels,
including zero, have complements. In this case for any pair of
labels $u,\bar{u}\in\rho_{\nu(p,r_0)}$, exactly one of these
labels is almost deterministic and, in particular, belongs to
$U_{\rm fin}$, and the other label marks a formula
$\theta_{p,\cdot,r_0}(a,y)$ with infinitely many solutions, where
$\models p(a)$, and belongs to $U_{\rm cofin}$. Hence, the
condition (c) holds.

If $r_0$ is non-isolated, then all non-zero labels, linking
realizations of $r_0$ are positive, since having a non-positive
non-zero label $u$, linking realizations of $r_0$ we have the
non-symmetric relation ${\rm SI}_{r_0}$ and as $r_0$ is
non-isolated there are infinitely many solutions for the formula
$\theta_u(x,a)$, where $\models r_0(a)$. This contradicts the
strong minimality of theory $T$. By Proposition 1.1 \cite{Su122},
non-zero labels linking realizations of types in
$R\setminus\{r_0\}$, are positive, and labels, linking
realizations of $r_0$ with realizations of types in
$R\setminus\{r_0\}$, are negative. In this case, since for
non-principal type there are only relative complements, if a label
$u$ belongs to $\rho_{\nu(p,r_0)}$, then $u$ is positive or zero,
almost deterministic and does not have a complement. Moreover,
$p=r_0$ since realizations of principal types cannot semi-isolate
realizations of non-principal type $r_0$. If the set
$\rho_{\nu(r_0)}$ is finite, then any label in $U_{\rm fin}$
belongs to $U_{\rm cofin}$ and vice versa, i.~e., $U_{\rm
fin}=U_{\rm cofin}$, and if $\rho_{\nu(r_0)}$ is infinite, then
all labels are almost deterministic and $U_{\rm
cofin}=\varnothing$. Thus, the condition (d) holds. The condition
(e) is implied by previous items.~$\Box$

\medskip
If $\mathfrak{M}$ is a ${\rm POSTC}_\mathcal{R}$-structure and
there is a theory $T$ with a family $R=S^1(\varnothing)$ and a
regular family $\nu(R)$ of labelling functions for semi-isolating
formulas such that $\mathfrak{M}_{\nu(R)}=\mathfrak{M}$, then we
say that $\mathfrak{M}$ {\em is represented}\index{${\rm
POSTC}_\mathcal{R}$-structure!represented} by $T$ and also say
that $T$ {\em represents}\index{Theory!repsesenting ${\rm
POSTC}_\mathcal{R}$-structure} the ${\rm
POSTC}_\mathcal{R}$-structure $\mathfrak{M}$. If all types of $R$
are realized in a model $\mathcal{N}$ of $T$, then we say that
$\mathfrak{M}$ {\em is represented}\index{${\rm
POSTC}_\mathcal{R}$-structure!represented} by $\mathcal{N}$.

Note that the syntactic representability of ${\rm
POSTC}_\mathcal{R}$-structure $\mathfrak{M}$ (by a theory) is
equivalent to the semantic representability of  $\mathfrak{M}$ (by
a model).

Notice also that there is a representation $T$ for the ${\rm
POSTC}_\mathcal{R}$-structure $\mathfrak{M}$ such that a label $u$
is almost deterministic if and only if $u$ does not force an
infinite set of solutions.

\medskip
{\bf Theorem 2.} {\em Let $\mathfrak{M}$ be a ${\rm
POSTC}_\mathcal{R}$-structure satisfying the following conditions:

{\rm (a)} $\mathfrak{M}$ consists of labels belonging to $U_{\rm
fin}\cup U_{\rm cofin}$;

{\rm (b)} there is an element $r_0\in \mathcal{R}$ such that any
set $\rho_{\mu(p,q)}$ is finite, where $p,q\in \mathcal{R}$ and
$q\ne r_0$, and all labels $u\in\rho_{\mu(p,q)}$ are almost
deterministic {\rm (}in some representation $\mathcal{N}$ of
$\mathfrak{M}${\rm )} and belong to $U_{\rm fin}$;

{\rm (c)} if $\mathcal{R}$ is finite then all non-zero labels are
positive and all labels $u$, including zero, have complements
$\bar{u}$, and for any pair of labels
$u,\bar{u}\in\rho_{\mu(p,r_0)}$, exactly one of them is almost
deterministic and, in particular, belongs to $U_{\rm fin}$, and
the other label marks a formula $\theta_{p,\cdot,r_0}(a,y)$ {\rm
(}for $mathcal{N}${\rm )} with infinitely many solutions, where
$\models p(a)$, and belongs to $U_{\rm cofin}$;

{\rm (d)} if $\mathcal{R}$ is infinite then all non-zero labels,
linking $r_0$ or elements of $R\setminus\{r_0\}$, are positive,
and labels, linking $r_0$ with elements in $R\setminus\{r_0\}$,
are negative; in this case, if a label $u$ belongs to
$\rho_{\mu(p,r_0)}$, then $u$ is positive or zero, almost
deterministic, does not have complements and $p=r_0$, moreover,
$U_{\rm fin}=U_{\rm cofin}$ if $\rho_{\mu(r_0)}$ is finite, and
$U_{\rm cofin}=\varnothing$ if $\rho_{\mu(r_0)}$ is infinite;

{\rm (e)} only labels in $\rho_{\mu(p,r_0)}$ and with
$|\mathcal{R}|<\omega$ can force infinity.

Then there is a strongly minimal theory $T$ representing the ${\rm
POSTC}_\mathcal{R}$-structure $\mathfrak{M}$ and having unique
$1$-type $r_0$ with infinitely many realizations.}

\medskip
{\em Proof.} Consider the construction for the proof of Theorems
9.1 \cite{ShS} and 8.1 \cite{Su122}. We identify $\mathcal{R}$
with the set of $1$-types for the required theory $T$. Now we add
to the types describing links between elements with respect to
binary relations $Q_u$, $u\in U$, an information for the
cardinality of sets of solutions for formulas
$\theta_{p,u,q}(a,y)$, where $\models p(a)$. The generic
construction for the class ${\bf T}_0$ of types guarantees that
the generic theory of the language $\{Q_u\mid u\in U\}$ is
strongly minimal and represents the ${\rm
POSTC}_\mathcal{R}$-structure $\mathfrak{M}$.~$\Box$


\bigskip
\begin{thebibliography}{1}
\bibitem{ShS}{\em Shulepov I.~V.} Algebras of distributions
for binary isolating formulas of a complete theory~/
I.~V.~Shulepov, S.~V.~Sudoplatov // arXiv:1205.3473v1
[math.LO].~--- 2012.~--- 41 p.
\bibitem{Su122}{\em Sudoplatov S.~V.} Algebras of distributions
for binary semi-isolating formulas of a complete theory~/
S.~V.~Sudoplatov // arXiv:1210.4049v1 [math.LO].~--- 2012.~---
35~p.
\bibitem{BaLa}{\em Baldwin J.~T.} On strongly minimal
sets~/ J.~T.~Baldwin, A.~H.~Lachlan // J.~Symbolic Logic.~---
1971.~--- Vol.~36, No.~1.~--- P. 79--96.
\bibitem{SuGP} {\em Sudoplatov S.~V.} Group polygonometries~/ S.~V.~Sudoplatov.~---
Novosibirsk~: NSTU, 2011, 2013.~--- 302~p. [in~Russian]
\end{thebibliography}
\end{document}